\documentclass[12pt]{article}
\usepackage{amsmath}
\usepackage{CJK}
\usepackage{amsfonts}
\usepackage{mathrsfs}
\usepackage{amssymb}
\usepackage[usenames]{color}

\def\g{\mathfrak{g}}
\def\b{\mathfrak{b}}
\def\a{\alpha}
\def\l{\lambda}
\def\m{\mu}
\def\D{\Delta}

\def\p{\partial}

\def\Vir{{\rm Vir}}
\def\sc{\scriptstyle}
\def\ssc{\scriptscriptstyle}
\def\dis{\displaystyle}
\def\cl{\centerline}

\def\vs{\vspace*}

\def\ptl{\partial}

\def\Z{\mathbb{Z}{\ssc\,}}

\def\C{\mathbb{C}{\ssc\,}}

\numberwithin{equation}{section}
\newtheorem{theo}{Theorem}[section]
\newtheorem{defi}[theo]{Definition}

\newtheorem{lemm}[theo]{Lemma}

\def\GG{{\cal G}}
\def\B{{\cal B}}
\def\Vir{{\rm Vir}}

\begin{document}
\begin{CJK*}{GBK}{song}

\begin{center}{\large \bf
Another class of simple graded Lie conformal algebras \\ that cannot be embedded into general Lie conformal algebras}
\footnote {Supported by NSF grant no.~11971350 of China.
}
\end{center}
\vs{6pt}

\cl{Yucai Su,  \ \ Xiaoqing Yue}
\cl{\footnotesize School of Mathematical Sciences, Tongji University, Shanghai 200092, China}

\cl{\footnotesize E-mails: ycsu@tongji.edu.cn, xiaoqingyue@tongji.edu.cn}
\vs{12pt} \par

{\small \parskip .005 truein \baselineskip 3pt \lineskip 3pt

\noindent{{\bf Abstract.} {In a previous paper by the authors, we obtain the first example
of a finitely freely generated simple $\mathbb Z$-graded Lie conformal algebra of linear
growth that cannot be embedded into any general Lie conformal algebra. In this paper,
we obtain, as a byproduct, another class of such
Lie conformal algebras by classifying $\mathbb Z$-graded simple Lie conformal algebras ${\cal G}=\oplus_{i=-1}^\infty{\cal G}_i$ satisfying
the following,\parskip-3pt
\begin{itemize}\parskip-3pt
\item[(1)]${\cal G}_0\cong{\rm Vir}$, the Virasoro conformal algebra;
\item[(2)]Each ${\cal G}_i$ for $i\ge-1$ is a ${\rm Vir}$-module of rank one.
\end{itemize}
These algebras include some Lie conformal algebras of Block type. }
 \vs{6pt}

}

}
\parskip .001 truein\baselineskip 8pt \lineskip 8pt

\vs{12pt}
\par

\cl{\bf\S1. \
Introduction}\setcounter{section}{1}\setcounter{equation}{0}\setcounter{theo}{0}

Conformal algebras, first introduced in [\ref{K1}], appear naturally in the context of formal distribution Lie algebras
and 
play important roles in quantum
field theory and conformal field theory (e.g., [\ref{BPZ}, \ref{CKW}, \ref{K3}]). They also turn out to be
 effective 
tools in the study of infinite-dimensional
Lie or associative algebras satisfying the
locality property, and their representations [\ref{K2}].

In recent years, the structure theory, representation theory and cohomology theory of Lie conformal
algebras have been extensively studied 
(e.g., [\ref{BKV}, \ref{BKLR}--\ref{KP2}, \ref{SXY}, \ref{Z1}]).
In particular,
simple finite Lie conformal algebras were classified in
[\ref{DK}], which turn out to be 
isomorphic either to
the Virasoro conformal algebra or  the
current Lie conformal algebra Cur$\,\g$ associated to a simple
finite-dimensional Lie algebra $\g$.
Finite irreducible conformal modules over the
Virasoro conformal algebra were determined
in [\ref{CK}] and of their extensions in
[\ref{CKW}]. The  cohomology theory of conformal
algebras with coefficients in an arbitrary module has been developed
in [\ref{BKV}, \ref{DSK}]. However, the theory of simple infinite Lie conformal algebras is
 far from being well developed, it is more complicated than the theory of Lie or associative algebras (e.g., [\ref{KP1}, \ref{KP2}]).

In order to better understand the theory of simple infinite Lie conformal
algebras, it is very natural to first study some important examples. It
is well-known that the general Lie conformal algebra $gc_N$ (which is a simple infinite Lie conformal
algebra) plays the same important role in the theory of Lie conformal algebras as the general Lie
algebra $gl_N$ does in the theory of Lie algebras. Thus
the study of Lie conformal algebras related to the general Lie conformal algebra $gc_N$ has drawn lots of attention in literature
(e.g., [\ref{BKL1}, \ref{BKL2}, \ref{GXY}, \ref{S1}--\ref{WCY}]).
In particular, in \cite{SY}, we study filtered Lie conformal algebras whose associated graded algebras are isomorphic to that of the general Lie
conformal algebra $gc_1$,  and as a byproduct we obtain the first example of a
finitely freely generated simple $\mathbb Z$-graded Lie conformal algebra of linear
growth that cannot be embedded into a general Lie conformal algebra $gc_N$ for any $N$, namely, the
Lie conformal algebra ${\rm gr\,}gc_1$ (the
associated graded conformal algebra of $gc_1$, which is also called a Lie conformal algebra of Block type), see \cite[Theorem 1.1]{SY}.
Motivated by the facts that
a simple Lie conformal algebra of rank one is isomorphic to the Virasoro conformal algebra $\Vir$
and that a finite simple $\Vir$-module is of rank one
[\ref{CK}--\ref{DK}], in this paper, we study
$\mathbb Z$-graded Lie conformal algebras ${\cal G}=\oplus_{i=-1}^\infty{\cal G}_i$ satisfying the following reasonable conditions\parskip-3pt
\begin{itemize}\parskip-3pt
\item[(C1)]${\cal G}_0\cong{\rm Vir}$, the Virasoro conformal algebra;
\item[(C2)]Each ${\cal G}_i$ for $i\ge-1$ is a ${\rm Vir}$-module of rank one;
\item[(C3)]${\cal G}$ is simple.
\end{itemize}
To state the main result, we first give the following definitions.
\begin{defi}\label{def-1}\rm
Let $\a\in\C$,\,$s=1,2$. Denote by $B(s,\a)$ the Lie conformal algebra with $\C[\partial]$-basis $\{G_i\,|\,i\in \Z_{\geq-1}\}$ and the following
$\l$-brackets, 
\begin{eqnarray}
\label{th11}
B(1,\a):&\!\!&
[G{_{-1}}_{\l}G_{-1}]=0,\, \  [G{_{-1}}_{\l}G_0]=(\a-\p) G_{-1},\, \
[G{_{-1}}_{\l}G_j]=(j+1) G_{j-1},\,  \ j\geq 1,\nonumber\\
&\!\!&
[G{_i}_\l G_j]=\Big((j-i)\a+(i+j+2)\l+(i+1)\p\Big)G_{i+j},\ \ i,j\in\Z_+,\nonumber\\
\label{th12}
B(2,\a):&\!\!&
[G{_i}_\l G_j]=\Big((j-i)\a+(i+j+2)\l+(i+1)\p\Big)G_{i+j},\ \ i,j\in\Z_{\geq-1}.
\end{eqnarray}
\end{defi}

\begin{defi}\label{def-2}\rm
Let $\a\in\C$. Denote by $\B(\a)=\oplus_{i=-1}^\infty \B_i$ the $\mathbb Z$-graded simple Lie conformal algebra with the $\l$-brackets $[\b{_0}_\l \b_i]=\big(i\a+(i+2)\l+\p\big)\b_{i}$ for $i> 0$, which satisfies
\begin{eqnarray}&\!\!\!\!\!\!\!\!\!\!\!\!\!\!\!\!\!\!\!\!\!\!\!\!&
{\rm(i)\ }{\cal B}_0\cong{\rm Vir},\ \ \ \ \mbox{the Virasoro conformal algebra,}\nonumber\\
&\!\!\!\!\!\!\!\!\!\!\!\!\!\!\!\!\!\!\!\!\!\!\!\!&
{\rm(ii)\ }{\cal B}_{-1} \mbox{\ is a Vir-module of rank one,}\nonumber\\
&\!\!\!\!\!\!\!\!\!\!\!\!\!\!\!\!\!\!\!\!\!\!\!\!&
{\rm(iii)\ }\mbox{Each}\  {\cal B}_i \mbox{\ for\ } i> 0\mbox{\ is a Vir-module of finite rank,}\nonumber
\end{eqnarray}
where $\b_0$ is the $\C[\ptl]$-generator of $\B_0$ and $\b_i$ is any one of $\C[\ptl]$-generators of $\B_i$ for $i> 0$.
\end{defi}

The main result of the present paper is the following.
\begin{theo}\label{main-theo}
\begin{itemize}
\item[\rm(1)]The Lie conformal algebra $B(s,\a)$ is simple for any $\a\in\C$ and $s=1,2$.
\item[\rm(2)]For $\a_1,\a_2\in\C,\,s_1,s_2\in\{1,2\}$, $B(s_1,\a_1)\cong B(s_2,\a_2)$ if and only if $(s_1,\a_1)=(s_2,\a_2)$.
\item[\rm(3)]Let $\GG=\oplus_{i=-1}^\infty\GG_i$ be a simple Lie conformal algebra satisfying conditions {\rm(C1)} and {\rm(C2)}. Then $\GG\cong B(s,\a)$ for some $\a\in\C$ and $s=1,2$.
\item[\rm(4)]For any $\a\in\C$, the Lie conformal algebra $\B(\a)$ does not have a nontrivial representation on any finite $\C[\partial]$-module.
In particular, $\B(\a)$ is a finitely freely generated simple Lie conformal algebra of linear
growth that cannot be embedded into $gc_N$ for any $N$.
\end{itemize}
\end{theo}
%
%
%
%
%

Therefore, Theorem \ref{main-theo}\,(4) provides another class $\B(\a)$ of
finitely freely generated simple $\mathbb Z$-graded Lie conformal algebras of linear
growth that cannot be embedded into a general Lie conformal algebra $gc_N$ for any $N$.

The paper is organized as follows. In section 2, we briefly recall some definitions
and preliminary results.  In section 3, we first study the structure of the Lie conformal algebra $B(s,\a)$, then we give the proof of Theorem \ref{main-theo}\,(4). In order to classify $\mathbb Z$-graded simple Lie conformal algebras ${\cal G}$, some technical lemmas were given in section 4. Then in section 5, we use these technical lemmas
to determine all
simple Lie conformal algebras satisfying conditions (C1), (C2), and complete the proof of Theorem \ref{main-theo}.

Throughout the paper, we denote by $\C,\,\C^*,\, \Z,\, \Z_+,\, \Z_{\geq -1}$ the sets of complex numbers, nonzero complex numbers, integers, nonnegative integers  and integers greater than $-2$
respectively.
\vs{12pt}

\cl{\bf\S2. \ Definitions and preliminary results }\setcounter{section}{2}\setcounter{equation}{0} \setcounter{theo}{0}

In this section, we summarize some basic definitions and results concerning Lie conformal algebras. More details can be found in [\ref{BKV}, \ref{CK}, \ref{K1}].

\begin{defi}\rm A
{\it Lie conformal algebra } is a $\C[\p]$-module $A$ with a
$\l$-bracket $[\cdot_{\l}\cdot]$ which defines a $\C$-bilinear map $A\times A
\rightarrow A[\l]$, where $A[\l]=\C[\l]\otimes A$ is the space of
polynomials of $\l$ with coefficients in $A$, such that for \vs{-5pt}$x,y,z\in A$,
\begin{eqnarray}&\!\!\!\!\!\!\!\!\!\!\!\!\!\!\!\!\!\!\!\!\!\!\!\!&\label{J}
[\p x_{\l}y]=-\l [x_\l y],\ \ \ \  [x_\l \p y]=(\p+\l)[x_\l y] \ \ \   \mbox{(conformal  sesquilinearity)},\\
&\!\!\!\!\!\!\!\!\!\!\!\!\!\!\!\!\!\!\!\!\!\!\!\!&\label{J_a_b_}
[x_{\l}y]=-[y_{-\l-\p} x] \ \ \   \mbox{(skew-symmetry)},\\
&\!\!\!\!\!\!\!\!\!\!\!\!\!\!\!\!\!\!\!\!\!\!\!\!&\label{J_a_b_c}
[x_{\l}[y_\m z]]=[[x_\l y]_{\l+\m} z]+[y_\m[x_\l z]] \ \ \
\mbox{(Jacobi identity)}.
\end{eqnarray}

A subset $S\subset A$ is called a {\it generating set} if $S$ generates $A$ as a $\C[\p]$-module. If there exists a finite generating set, then $A$ is called {\it finite}. Otherwise, it is called {\it infinite}.
\end{defi}

For a given Lie conformal algebra $A$, from \cite{K1}, we know that there is an important Lie algebra associated to it. For each $j\in \mathbb{Z}_{+}$, regarding $[a_\lambda b]\in \C[\lambda]\otimes A$ as a formal polynomial in $\lambda$, we can define the \emph{$j$th product} $a_{(j)}b$ by the coefficient of $\lambda^j$ in $[a_\lambda b]$, i.e. $a_{(j)}b$ for all $a,b\in A$ as follows:
\begin{equation}\label{Tl3.1}
[a_\l b]=\sum_{j\in \mathbb{Z}_+}(a_{(j)}b)\frac{\l^{j}}{j!}.
\end{equation}
Now we can give the definition of this Lie algebra.
\begin{defi}\rm An
{\it annihilation algebra } of a Lie conformal algebra $A$ is a Lie algebra  with $\C$-basis $\{a_{(n)} \mid a\in A,\, n\in \Z_+\}$ and relations $$[a_{(m)}, b_{(n)}]=\sum_{j\in \mathbb{Z}_+}{m \choose j}(a_{(j)}b)_{(m+n-j)},\ \ \ \ \  \p(a_{(n)})=-na_{(n-1)}.$$
\end{defi}

The {\it Virasoro conformal algebra} $\Vir$ is the simplest nontrivial Lie conformal algebra. It is a free $\C[\p]$-module of rank one with generator $L$ and can be defined by
\begin{equation}\label{Vir}
\Vir=\C[\ptl]L:\ \ \ \ [L_\l L]=(\ptl+2\l)L.
\end{equation}
It is known that any simple Lie conformal algebra of free rank one over $\C[\p]$ is isomorphic to $\Vir$ [\ref{DK}].

The {\it general Lie conformal algebra} $gc_N$ can be defined as the
infinite rank
$\C[\p]$-module  $\C[\ptl,x]\otimes gl_N$,
 with the $\l$-bracket
\begin{equation}\label{gc-N}
[{f(\ptl,x)A}_{\,\l\,} g(\ptl,x)B]{\sc\!}={\sc\!}
f(-\l,x{\sc\!}+{\sc\!}\ptl{\sc\!}+{\sc\!}\l)
g(\ptl{\sc\!}+{\sc\!}\l,x)AB{\sc\!}-{\sc\!}f(-\l,x)g(\ptl{\sc\!}
+{\sc\!}\l,x{\sc\!}-{\sc\!}\l)BA,
\end{equation}
for $f(\ptl,x),g(\ptl,x)\in\C[\ptl,x],$ $A,B\in gl_N$, where $gl_N$ is the space of $N\times N$ matrices, and we have identified $f(\ptl,x)\otimes A$ with $f(\ptl,x)A$.
If we set $J_A^n=x^nA$, \vs{-10pt}then
\begin{equation*}
[{J_{A}^{m}}_{\l} J_B^{n}]=\mbox{$\sum\limits_{s=0}^{m}$}\binom{m}{s}
(\l+\p)^{s}J_{AB}^{m+n-s}-\mbox{$\sum\limits_{s=0}^{n}$} \binom{n}{s}
(-\l)^{s}J_{BA}^{m+n-s}\vs{-5pt},
\end{equation*}
for $m,n\in \Z_+,$ $A,B\in gl_N,$ where
$\big(^m_{\,s}\big)=m(m-1)\cdot\cdot\cdot(m-s+1)/{s}!$ if $s\geq 0$
and $\big(^m_{\,s}\big)=0$ otherwise, is the binomial coefficient.

\begin{defi}\rm A {\it module over a Lie conformal algebra} $A$ is a $\C[\p]$-module $M$ with a $\l$-action $\cdot_{\l}\cdot: A\times M\rightarrow M[[\l]]$, where $M[[\l]]$ is the set of formal power series of $\l$ with coefficients in $M$, such that for $x,y \in A,\ v\in M$,
\begin{eqnarray}
&&x_{\l}(y_{\mu}v)-y_{\mu}(x_{\l}v)=[x_{\l}y]_{\l+\mu}v,\\
&&(\p x)_{\l}v=-\l x_{\l}v,\ \ x_{\l}(\p v)=(\p+\l)x_{\l}v.
\end{eqnarray}
An $A$-module $M$ is called {\it conformal} if $x_{\l}v\in M[\l]$ for $x\in A,\ v\in M$ and {\it finite} if $M$ is finitely generated over $\C[\p]$.
\end{defi}

According to [\ref{CK}], we know that all free nontrivial conformal $\Vir$-modules of rank one over $\C[\p]$ are $M_{\Delta,\alpha}$ for $\Delta,\alpha\in \C$, where
\begin{equation}\label{M-a-b}
M_{\Delta,\alpha}=\C[\ptl]v:\ \ \ \ L_\l v=(\alpha+\ptl+\Delta\l)v.
\end{equation}
The module $M_{\Delta,\alpha}$ is irreducible if and only if $\Delta\ne0$.
\vs{12pt}

\cl{\bf\S3. \  Graded Lie conformal algebras $B(s,\a)$}\setcounter{section}{3}\setcounter{equation}{0}\setcounter{theo}{0} \vs{5pt}

It is straightforward to verify that \eqref{th11} indeed defines Lie conformal algebras $B(1,\a)$ and $B(2,\a)$. Furthermore,
the annihilation algebra of $B(2,\a)$ is the Lie algebra ${\cal A}=\mbox{span}_{\C}\{ G_{i,m}\,|\,i,m\in\Z_{\geq-1}\}$ with Lie brackets $$[G_{i,m}, G_{j,n}]=(j-i)\a G_{i+j,m+n+1}+\big((j+1)(m+1)-(n+1)(i+1)\big)G_{i+j,m+n}.$$ When $\a=0$, this Lie algebra has close relation to the Block-type Lie algebras studied in [\ref{GXY}, \ref{SXY}].

Now we study the structure of the Lie conformal algebra $B(s,\a)$ for $\a\in\C$ and $s=1,2$. First we need some definitions.
For any $x\in B(s,\a)$, we define the operator $({\rm ad\,}x)_\l:B(s,\a)\to B(s,\a)[\l]$ such that $({\rm ad\,}x)_\l(y)=[x_\l y]$ for any $y\in B(s,\a)$.
An element $x\in B(s,\a)$ is {\it locally nilpotent} if for any $y\in B(s,\a)$, there exists $1\le n\in\Z_+$ such that $({\rm ad\,}x)_\l^n(y)=0$. We have
\begin{lemm}\label{nilpo}The set
 of locally nilpotent elements of $B(s,\a)$ is equal to $\C[\p]G_{-1}$. \end{lemm}\noindent{\it Proof.~}~Denote by ${\cal N}$ the set of locally nilpotent elements of $B(s,\a)$. First by \eqref{th11} and conformal  sesquilinearity, we have the following for
$a_0(\p),c_j(\p)\in\C[\p]$,
\begin{eqnarray}\label{Thefff}
[a_0(\p)G{_{-1}}_{\l}c_j(\p)G_j]=\left\{\begin{array}{lll}
(\a-\p)a_0(-\l)c_j(\p+\l)G_{-1}&\mbox{if \ }s=1,\ j=0,\\[4pt]
(j+1)a_0(-\l)c_j(\p+\l) G_{j-1}&\mbox{if \ }s=1,\ 1\leq j\in\Z,\\[4pt]
(j+1)(\a+\l)a_0(-\l)c_j(\p+\l) G_{j-1}&\mbox{if \ }s=2,\ 0\leq j\in\Z.
\end{array}\right.
\end{eqnarray}
From this and using conformal sesquilinearity, we immediately obtain that $\C[\p]G_{-1}\subset{\cal N}$.
Now let $x=\sum_{i=-1}^{\infty}b_i(\p)G_i\in {\cal N}$, suppose $max\{i\,|\,b_i(\p)\neq 0\}=i_0$. If $i_0\geq 0$, by \eqref{th11},
we have $[G{_{i_0}}_\l G_j]=\big((j-{i_0})\a+({i_0}+j+2)\l+({i_0}+1)\p\big)G_{{i_0}+j}$ for $j\in\Z_+$.
Then applying $({\rm ad\,}x)_\l^{n}$ to $G_j$, we can obtain that the coefficient of $G_{ni_0+j}$ in the expression of $({\rm ad\,}x)_\l^n(G_j)$ is nonzero for any $1\le n\in\Z_+$ and $j\in\Z_+$. This is a contradiction with $x\in {\cal N}$. Therefore we must have $i_0=-1$, then the lemma follows.
\hfill$\Box$\vs{7pt}

\noindent{\it Proof of Theorem \ref{main-theo}\,(1) and (2).~}~(1)~Let $J$ be a nonzero ideal of $B(s,\a)$ for some $\a\in\C$ and $s=1,2$. Then there exists at least one nonzero element $x=\sum_{j=-1}^{m}b_j(\p)G_j\in J$ for some $b_j(\p)\in\C[\partial]$, where $m\in\Z_+$  such that $b_m(\p)\neq 0$.
We claim that $a_0(\p)G_{-1}\in J$ for some nonzero $a_0(\p)\in\C[\partial]$. If $m=-1$, we immediately have the claim. Otherwise, we can apply the operator $({\rm ad\,}G_{-1})_\l^{m+1}$ to $x$, we have
the following for
$b_m(\p)\in\C[\p]$,
\begin{eqnarray}
J\ni({\rm ad\,}G_{-1})_\l^{m+1}(x)=\left\{\begin{array}{ll}
(m+1)!(\a-\p)b_m(\p+\l)G_{-1}&\mbox{if \ }s=1,\\[4pt]
(m+1)!(\a+\l)^{m}b_m(\p+\l) G_{-1}&\mbox{if \ }s=2.
\end{array}\right.
\end{eqnarray}
Then we inductively deduce from \eqref{Thefff} that all $G_{j}\in J$ for $j\in \Z_{\geq-1}$, i.e., $J=B(s,\a)$. Therefore, $B(s,\a)$ is simple.

(2)~It is obvious that the sufficient condition holds. We only need to prove the necessary condition. For $\a_1,\a_2\in\C,\,s_1,s_2\in\{1,2\}$, we suppose $B(s_1,\a_1)\cong B(s_2,\a_2)$, $\{G_i\,|\,i\in\Z_{\geq-1}\}$ and $\{G'_i\,|\, i\in\Z_{\geq-1}\}$ are the $\C[\partial]$-bases of $B(s_1,\a_1)$ and $B(s_2,\a_2)$ respectively. By (\ref{th11}), we can immediately conclude that $s_1=s_2$.

First suppose $s_1=1$. Let $\varphi : B(1,\a_1)\longrightarrow B(1,\a_2)$ be an isomorphism. By Lemma \ref{nilpo}, we can assume $\varphi(G_{-1})=a(\p)G_{-1}'$ and $\varphi(G_{0})=\sum_{i=-1}^{\infty}b_i(\p)G_{i}'$ for some $a(\p),\ b_i(\p)\in\C[\p]$ with $i\in\Z_{\geq-1}$. Applying the isomorphism $\varphi$ to the both sides of $[G{_{-1}}_{\l}G_{0}]=(\a_1-\p)G_{-1}$,
comparing the coefficients of $G_{i-1}'$ for $1\leq i\in\Z$ and $G_{-1}'$ respectively,
we can deduce that $b_i(\p)=0$ for $1\leq i\in\Z$ and
\begin{equation}\label{iso1}
(\a_2-\p)a(-\l)b_{0}(\l+\p)=(\a_1-\p)a(\p).
\end{equation}
By \eqref{th11}, we have $[G{_{0}}_{\l}G_{0}]=(2\l+\p)G_{0}$. Applying the isomorphism $\varphi$ to this equation, then comparing the coefficients of $G_{0}'$, we can obtain that
\begin{equation}\label{iso2}
b_{0}(-\l)b_{0}(\l+\p)=b_{0}(\p).
\end{equation}
Comparing the degrees of $\l$ in \eqref{iso2} and using the fact that $\varphi$ is an isomorphism, we have $b_{0}(\p)=1$.
Then by \eqref{iso1}, we can conclude that $\a_1=\a_2$.

Now suppose $s_1=2$. Similarly, if $B(2,\a_1)\cong B(2,\a_2)$, we can also obtain that $\a_1=\a_2$. Therefore, if $\varphi : B(s,\a_1)\longrightarrow B(s,\a_2)$ is an isomorphism, then $(s_1, \a_1)=(s_2, \a_2)$.
\hfill$\Box$

In order to prove Theorem \ref{main-theo} (4), we need some preparations.  Assume $V$ is
a finitely freely $\C[\ptl]$-generated nontrivial $\B(\a)$-module. Regarding $V$ as a module over $\Vir$, by \cite[Theorem 3.2(1)]{CK}, we can choose a composition \vs{-5pt}series,
$$
V=V_{N}\supset V_{N-1}\supset\cdots\supset V_{1}\supset V_{0}=0\vs{-5pt},
$$
such that for each $i=1,2,...,N$, the composition factor $\overline{V}{}_{i}=V_{i}/V_{i-1}$ is either a rank one free module $M_{\Delta_i,\beta_i}$ with $\D_i\ne0$, or else a $1$-dimensional trivial module $\C_{\beta_i}$ with trivial $\l$-action and with $\ptl$ acting as the  scalar $\beta_i$. Denote by $\bar v_i$ a $\C[\ptl]$-generator of $\overline{V}{}_{i}$ and $v_i\in V_i$ the preimage of $\bar v_i$.
Then $\{v_i\,|\,1\le i\le N\}$ is a $\C[\ptl]$-generating set of $V$, such that the $\l$-action of $\b_0$ on $v_i$ is a $\C[\l,\ptl]$-combination of $v_1,...,v_i$.
\begin{lemm}\label{J-1-l-ac}For all $i\gg0$,  the $\l$-action of $\b_i$ on $v_1$ is trivial, namely, ${\b_i}_{\l}v_1=0$.
\end{lemm}
\noindent{\it Proof.~}~Assume  $i\gg0$ is fixed and suppose ${\b_i}_{\l}v_1\ne0$, and let $k_i\ge1$ be the largest integer such that ${\b_i}_{\l}v_1\not\subset V_{k_i-1}[\l]$. We consider the following possibilities.\vs{5pt}

\noindent{\bf Case 1:} $V_1=M_{\D_1,\beta_1},\,\overline V_{k_i}=M_{\D_{k_i},\beta_{k_i}}$.\vs{5pt}

We can 
write
\begin{equation}\label{SSS}
{\b_i}_\l v_1\equiv p_i(\l,\p)v_{k_i}\ ({\rm mod\,}V_{{k_i}-1}[\l])\mbox{ \ for some \ }p_i(\l,\p)\in\C[\l,\p]\vs{-5pt}.
\end{equation}
Applying the operator ${\b_0}_\mu$ to \eqref{SSS}, we \vs{-5pt}obtain
\begin{equation}\label{SSS1}
(\beta_{k_i}\!+\!\ptl\!+\!\D_{k_i}\mu)p_i(\l,\mu\!+\!\p)
\!=\!\big(i\a+(1\!+\!i)\mu\!-\!\l\big)p_i(\l+\mu,\p)\!+\!(\beta_1\!+\!\l\!+\!\p\!+\!\D_1\mu)p_i(\l,\p)\vs{-5pt}.
\end{equation}
Letting $\p=0$, we have
\begin{equation}\label{==SSS1-PP}
p_i(\l,\mu)=\frac1{\beta_{k_i}+\D_{k_i}\mu}
\Big(\big(i\a+(1+i)\mu-\l\big)p_i(\l+\mu,0)+(\beta_1+\l+\D_1\mu)p_i(\l,0)\Big)\vs{-10pt}.
\end{equation}
Using this in \eqref{SSS1} with $\l=i\a+(1+i)\mu$ and $\p=-\beta_{k_i}-\D_{k_i}\mu$, we can deduce that
$$\begin{array}{ll}
(i+1)\big((\D_{k_i}{\sc}+1)\mu+\beta_{k_i}\big)p_i\big((i+1-\D_{k_i})\mu-\beta_{k_i}+i\a,0\big)\\[5pt]
\ \ \ \ =\big((i+1-\D_1\D_{k_i})\mu+\beta_1-\beta_{k_i}\D_1\big)p_i\big((i+1)\mu+i\a,0\big)\vs{-5pt}.
\end{array}$$
Suppose $p_i(\l,0)$ has degree $m_i$. Comparing the coefficients of $\mu^{m_i+1}$ in the above equation, we obtain (note that the following equation does not depend on the coefficients of $p_i(\l,\mu)$\vs{-5pt})
\begin{equation}\label{==SSS1-PP===}
(i\!+\!1)(\D_{k_i}\!+\!1)(i\!+\!1\!-\!\D_{k_i})^{m_i}
\!=\!(i\!+\!1\!-\!\D_1\D_{k_i})(i\!+\!1)^{m_i}\vs{-5pt}.
\end{equation}
When $i$ is sufficiently large, one can easily see that \eqref{==SSS1-PP===} cannot hold if $m_i>1$ (note that $\D_1,\D_{k_i}\ne0$, and there are only finitely many choices for what $\D_{k_i}$ can be, since $1\le k_i\le N$). Thus $m_i\le1$ if $i\gg0$.
 Then from \eqref{==SSS1-PP}, we obtain that $p_i(\l,\mu)$ is a polynomial of degree $\le1$. Thus suppose $p_i(\l,\mu)=a_{i,0}+a_{i,1}\l+a_{i,2}\mu$. Then by comparing the coefficients of $\l\m$, $\m\p$ and $\m$ respectively in \eqref{SSS1}, we immediately have $p_i(\l,\mu)=0$.\vs{5pt}

\noindent{\bf Case 2:} $V_1=\C_{\beta_1},\,\overline V_{k_i}=M_{\D_{k_i},\beta_{k_i}}$.\vs{5pt}

In this case, we can still assume \eqref{SSS}. Applying the operator ${\b_0}_\mu$ to \eqref{SSS}, we \vs{-5pt}obtain
\begin{equation}\label{SSS2}
p_i(\l,\mu+\p)(\beta_{k_i}+\ptl+\D_{k_i}\mu)=(i\a+(1+i)\mu-\l)p_i(\l+\mu,\p)\vs{-5pt}.
\end{equation}
Taking $\mu=\p=0$, we get $p_i(\l,0)=0$. Then letting $\p=0$, we obtain $p_i(\l,\mu)=0$.\vs{5pt}

\noindent{\bf Case 3:} $V_1=M_{\D_1,\beta_1},\,\overline V_{k_i}=\C_{\beta_{k_i}}$.\vs{5pt}

In this case, since   $\p$ acts on $\bar v_{k_i}$ as the scalar $\beta_{k_i}$, i.e., $\p v_{k_i}\equiv \beta_{k_i} v_{k_i}\,({\rm mod\,}V_{{k_i}-1}[\l])$, we can writ\vs{-5pt}e
\begin{equation}\label{SSS3}
{\b_i}_\l v_1\equiv p_i(\l)v_{k_i}\ ({\rm mod\,}V_{{k_i}-1}[\l])\mbox{ \ for some \ }p_i(\l)\in\C[\l].
\end{equation}
Applying the operator ${\b_0}_\mu$ to \eqref{SSS3}, we \vs{-5pt}obtain
\begin{equation}\label{SSS3-1}
0=\big(i\a+(1+i)\mu-\l\big)p_i(\l+\mu)+(\beta_1+\l+\p+\D_1\mu)p_i(\l)\vs{-5pt}.
\end{equation}
By comparing the coefficients of $\p$, we immediately get $p_i(\l)=0.$\vs{5pt}

\noindent{\bf Case 4:} $V_1=\C_{\beta_1},\,\overline V_{k_i}=\C_{\beta_{k_i}}$.\vs{5pt}

As above, we can still assume \eqref{SSS3}. Applying the operator ${\b_0}_\mu$ to \eqref{SSS3}, in this case we \vs{-5pt}obtain
\begin{equation}\label{SSS4-1}
0=\big(i\a+(1+i)\mu-\l\big)p_i(\l+\mu)\vs{-5pt}.
\end{equation}
It is obvious that $p_i(\l)=0$.
\hfill$\Box$\vs{5pt}

Finally we can give the proof of  Theorem \ref{main-theo} (4). By induction on $j\le N$, we obtain ${\b_i}_\l v_j=0$, i.e., the $\l$-action of $\b_i$ is trivial. From this, we immediately obtain that the $\l$-action of $\B(\a)$ on $V$ is trivial since
$\B(\a)$ is a simple Lie conformal algebra. Then Theorem \ref{main-theo} (4) follows.
\vs{10pt}

\cl{\bf\S4. \  Some technical lemmas}\setcounter{section}{4}\setcounter{equation}{0}\setcounter{theo}{0} \vs{5pt}
Let $\GG=\oplus_{i=-1}^\infty\GG_i$ be a simple Lie conformal algebra satisfying conditions {\rm(C1)} and {\rm(C2)}.
The main problem to be addressed in this paper is to
classify these $\mathbb Z$-graded Lie conformal algebras.
In order to solve the main problem, we need some preparations.
Since $\GG=\oplus_{i=-1}^{\infty}\GG_{i}$ satisfying conditions {\rm(C1)} and {\rm(C2)},
denote by $\g_i$ a $\C[\ptl]$-generator of $\GG_{i}$,
then $\{\g_i\,|\, i\in \Z_{\geq-1}\}$ is a $\C[\ptl]$-generating set of $\GG$.  By $(\ref {M-a-b})$, we can suppose
\begin{eqnarray}\label{G_-1_-1}
\!\!\!\!\!\!\!\!\!\!\!\!&\!\!\!\!\!\!\!\!\!\!\!\!\!\!\!&
[\g{_{-1}}_{\l}\g_{-1}]=0,
\\\!\!\!\!\!\!\!\!\!\!\!\!&\!\!\!\!\!\!\!\!\!\!\!\!\!\!\!&
[\g{_{0}}_{\l}\g_j]=(\alpha_j+\ptl+\Delta_j\l) \g_{j},\ \label{G_0_j}
\\\!\!\!\!\!\!\!\!\!\!\!\!&\!\!\!\!\!\!\!\!\!\!\!\!\!\!\!&
[\g{_i}_\l \g_j]=g_{i,j}(\l,\p)\g_{i+j},\label{G_i_j}
\end{eqnarray}
where $\alpha_j,\ \Delta_j\in\C$ for $j\in\Z_{\geq-1}$ and $g_{i,j}(\l,\p)\in\C[\l,\p]$ for $i,j\in\Z_{\geq-1}$ are polynomials of $\l$ and $\p$.
From $(\ref{G_i_j})$, we can see that $\g_{i+j}$ can be generated by $\g{_i}$ and $\g_j$ for $i,j\in\Z_{\geq-1}$. It is very natural to firstly consider $g_{i,j}(\l,\p)$ for $i=0$, $j\in\Z_{\geq-1}$ and $i=1$, $j=-1$. This is also the aim of this section. Based on this, in the next section we will determine all $g_{i,j}(\l,\p)$ for $i,j\in\Z_{\geq-1}$ in the proof of Theorem \ref{main-theo}\,(3). Since $\GG$ is simple, we also know that
\begin{equation}\label{G_-1_j}
[\g{_{-1}}_{\l}\g_{j}]\neq 0 \ \ \ \mbox{for}\ \ 1\leq j\in\Z.
\end{equation}

\begin{lemm}\label{lemm1}
In $(\ref{G_0_j})$, we have $\Delta_0=2$ and $\alpha_j=j\alpha_1$ for $ j\in\Z_{\geq-1}$; in particular, we get $\alpha_0=0$. Thus for $j\in\Z_{\geq-1}$ we can suppose
\begin{eqnarray}\label{g-0-0}
\!\!\!\!\!\!\!\!\!\!\!\!&\!\!\!\!\!\!\!\!\!\!\!\!\!\!\!&
g_{0,0}(\l,\p)=\ptl+2\l,
\\\!\!\!\!\!\!\!\!\!\!\!\!&\!\!\!\!\!\!\!\!\!\!\!\!\!\!\!&
g_{0,j}(\l,\p)=j\alpha_1+\ptl+\Delta_j\l.\ \label{g_0_j}
\end{eqnarray}
\end{lemm}

\noindent{\it Proof.~}~From ${\cal G}_0\cong{\rm Vir}$, by $(\ref{Vir})$, we can conclude that $\Delta_0=2$, $\alpha_0=0$ in $(\ref{G_0_j})$, thus $(\ref{g-0-0})$ holds.
Now applying the operator ${\g_0}_{\mu}$ to \eqref{G_i_j}, using the Jacobi identity $[\g{_0}_{\m}[\g{_{i}}_\l \g_{j}]]=[[\g{_0}_\m \g_{i}]_{\l+\m} \g_{j}]+[\g{_{i}}_\l[\g{_0}_\m \g_j]]$ and comparing the coefficients of $\g_{i+j}$ for $i,j\in\Z_{\geq-1},$ we obtain
\begin{eqnarray}
\label{aj}
\!\!\!\!\!\!\!\!\!\!\!\!&\!\!\!\!\!\!\!\!\!\!\!\!\!\!\!&\ \ \ \ \ \ \ \
(\a_{j+i}+\p+\D_{j+i}{\ssc\,}\m){\ssc\,}g_{i,j}(\l,\m+\p)-(\a_{j}+\p+\l+\D_{j}{\ssc\,}\m){\ssc\,}g_{i,j}(\l,\p)
\nonumber\\\!\!\!\!\!\!\!\!\!\!\!\!&\!\!\!\!\!\!\!\!\!\!\!\!\!\!\!&\ \ \ \ \ \ \ \
\ \ \ \ \ \ =\big(\a_{i}-\l+(\D_{i}-1){\ssc\,}\m\big){\ssc\,}g_{i,j}(\l+\m,\p).
\end{eqnarray}
Now taking $\m=0$ in \eqref{aj}, one can immediately get that $(\a_{j+i}-\a_{j}-\a_{i}){\ssc\,}g_{i,j}(\l,\p)=0$ for all $i,j\in\Z_{\geq-1}$. Then the lemma follows.
\hfill$\Box$\vspace{5pt}

In order to determine all the polynomials $g_{i,j}(\l,\p)$,  we would first like to deal with the case with $i=1$ and $j=-1$. There is no need to compute $g_{-1,1}(\l,\p)$ as it can be determined from
$g_{1,-1}(\l,\p)$ by skew-symmetry. Comparing the coefficients of $\g_{j+1}$ on both sides of the Jacobi identity
$[\g{_0}_{\l}[\g{_{1}}_\m \g_{j}]]=[[\g{_0}_\l \g_{1}]_{\l+\m} \g_{j}]+[\g{_{1}}_\m[\g{_0}_\l \g_j]],$  by
$(\ref{G_i_j})$ and Lemma \ref{lemm1}, we have
\begin{eqnarray}
\label{g1j}
\!\!\!\!\!\!\!\!&\!\!\!\!\!\!\!\!&
\big((j+1)\a_{1}+\p+\D_{j+1}{\ssc\,}\l\big){\ssc\,}g_{1,j}(\m,\l+\p)-(j\a_{1}+\p+\m+\D_{j}{\ssc\,}\l){\ssc\,}g_{1,j}(\m,\p)
\nonumber\\\!\!\!\!\!\!\!\!\!\!\!\!&\!\!\!\!\!\!\!\!\!\!\!\!\!\!\!&\ \ \ \ \ \ \ \
\ \ \ \ \ \ =\big(\a_{1}-\m+(\D_{1}-1){\ssc\,}\l\big){\ssc\,}g_{1,j}(\l+\m,\p).
\end{eqnarray}
Taking $\p=0$ and $j=-1$, by Lemma \ref{lemm1}, we obtain
\begin{eqnarray}\!\!\!\!\!\!\!\!\!\!\!\!&\!\!\!\!\!\!\!\!\!\!\!\!\!\!\!&
\label{g1 -1}\ \ \ \ \ g_{1,-1}(\m,\l)=\frac1{2\l}
\Big(\big(\a_{1}-\m+(\D_{1}-1){\ssc\,}\l\big){\ssc\,}g_{1,-1}(\l+\m,0)
\nonumber\\\!\!\!\!\!\!\!\!\!\!\!\!&\!\!\!\!\!\!\!\!\!\!\!\!\!\!\!&
\ \ \ \ \ \ \ \ \ \ \ \ \ \ \ \ \ \ \ \ \ \ \ -(\a_{1}-\m-\D_{-1}{\ssc\,}\l){\ssc\,}g_{1,-1}(\m,0)\Big).
\end{eqnarray}
Similarly, applying the operator ${\g_{-1}}_{\l}$ to $[\g{_1}_\m \g_{j}]$, using the Jacobi identity and comparing the coefficients of $\g_{j}$, by
$(\ref{G_i_j})$ and Lemma \ref{lemm1}, we obtain
\begin{eqnarray}
\label{g-1 1 j}
\!\!\!\!\!\!\!\!&\!\!\!\!\!\!\!\!&
g_{-1,j+1}(\l,\p){\ssc\,}g_{1,j}(\m,\l+\p)-g_{-1,j}(\l,\p+\m){\ssc\,}g_{1,j-1}(\m,\p)
\nonumber\\\!\!\!\!\!\!\!\!\!\!\!\!&\!\!\!\!\!\!\!\!\!\!\!\!\!\!\!&\ \ \ \ \ \ \ =\big(j\a_{1}+\p+\D_{j}{\ssc\,}(\l+\m)\big){\ssc\,}g_{-1,1}(\l,-\l-\m).
\end{eqnarray}
Setting $j=-1$ and replacing $\p$, $\l$ by $-\l$, $-\m-\p$ in $(\ref{g-1 1 j})$ respectively, by $(\ref{J_a_b_})$, $(\ref{G_-1_-1})$, $(\ref{G_i_j})$ and Lemma \ref{lemm1}, we can deduce that
\begin{equation}
\label{g1 -1 0}\ \ \ \ \ (\a_{1}-\m-\p) g_{1,-1}(\m,0)
=\big(\a_{1}-\m+(\D_{-1}-1{\ssc\,})\p\big){\ssc\,}g_{1,-1}(\m,\p).
\end{equation}
Using \eqref{g1 -1} in the above formula, we get
\begin{eqnarray}
\label{g1 -1 00}
\!\!\!\!\!\!\!\!&\!\!\!\!\!\!\!\!&
\big(\a_{1}-\m+(\D_{1}-1{\ssc\,})\p\big)\big(\a_{1}-\m+(\D_{-1}-1{\ssc\,})\p\big){\ssc\,}g_{1,-1}(\m+\p,0)
\nonumber\\\!\!\!\!\!\!\!\!\!\!\!\!&\!\!\!\!\!\!\!\!\!\!\!\!\!\!\!&\ \ \ \ \ \
=\Big(2\p(\a_{1}-\m-\p)+(\a_{1}-\m-\D_{-1}\p)\big(\a_{1}-\m+(\D_{-1}-1)\p\big)\Big){\ssc\,}g_{1,-1}(\m,0).
\end{eqnarray}
Letting $\m=0$ implies the following, 
\begin{eqnarray}
\label{g1 -1 01}
\!\!\!\!\!\!\!\!&\!\!\!\!\!\!\!\!&
\big(\a_{1}+(\D_{1}-1{\ssc\,})\p\big)\big(\a_{1}+(\D_{-1}-1{\ssc\,})\p\big){\ssc\,}g_{1,-1}(\p,0)
\nonumber\\\!\!\!\!\!\!\!\!\!\!\!\!&\!\!\!\!\!\!\!\!\!\!\!\!\!\!\!&\ \ \ \ \ \
=\Big(2\p(\a_{1}-\p)+(\a_{1}-\D_{-1}\p)\big(\a_{1}+(\D_{-1}-1)\p\big)\Big){\ssc\,}g_{1,-1}(0,0).
\end{eqnarray}
Taking $\p=-\m$ in $(\ref{g1 -1 00})$, and replacing $\m$ by $\p$, we obtain 
\begin{eqnarray}
\label{g1 -1 02}
\!\!\!\!\!\!\!\!&\!\!\!\!\!\!\!\!&
\Big(-2\p\a_{1}+(\a_{1}-\D_{-1}\p)\big(\a_{1}+(\D_{-1}-1)\p\big)\Big){\ssc\,}g_{1,-1}(\p,0)
\nonumber\\\!\!\!\!\!\!\!\!\!\!\!\!&\!\!\!\!\!\!\!\!\!\!\!\!\!\!\!&\ \ \ \ \ \
\ \ \ \ =(\a_{1}-\D_{1}{\ssc\,}\p)(\a_{1}-\D_{-1}{\ssc\,}\p){\ssc\,}g_{1,-1}(0,0).
\end{eqnarray}
Since $g_{1,-1}(\l,\p)$ is a polynomial of $\l$ and $\p$, we can write  $g_{1,-1}(\p,0)=\sum_{i=0}^{m}a_{1,-1}^{i}\p^{i}$ for some $a_{1,-1}^{i}\in\C$ for $0\leq i\leq m$.
We need to
consider whether or not $\a_1\ne0$.
First assume $\a_1\neq 0$. If
$a_{1,-1}^{0}=0$, then comparing the degrees of $\p$ on both sides of  $(\ref{g1 -1 02})$, we immediately have $g_{1,-1}(\p,0)=0$, which implies that  $g_{1,-1}(\l,\p)=0$ by $(\ref{g1 -1})$, a contradiction with $(\ref{G_-1_j})$. If $a_{1,-1}^{0}\ne0$, we have the following.

\begin{lemm}\label{lemm2}
Assume $\a_1\neq 0$ and $a_{1,-1}^{0}\neq 0$. Then $\Delta_{-1}=0,\ \Delta_1=3$ or $\Delta_{-1}=1,\ \Delta_1=3$. Furthermore, 
\begin{equation}\label{eq-1}
g_{1,-1}(\l,\p)=
\begin{cases}\displaystyle
a_{1,-1}^{0} \quad &\mbox{if \ \ $\Delta_{-1}=0,\ \Delta_1=3$,}\\[10pt]
\dis\frac{a_{1,-1}^{0}}{\a_1}\big(\a_1-\l-\p\big) \quad &\mbox{if \ \ $\Delta_{-1}=1,\ \Delta_1=3$.}
\end{cases}
\end{equation}
\end{lemm}
\noindent{\it Proof.~}~
Since $\a_1\neq 0$ and $g_{1,-1}(0,0)=a_{1,-1}^{0}\neq 0$,
comparing the degrees of $\p$ on both sides of $(\ref{g1 -1 02})$,
we get $a_{1,-1}^{i}=0$ for $2\leq i\leq m$. And
we need to consider the following possibilities.\vs{5pt}

\noindent{\bf Case 1:} $\Delta_{-1}=0$.\vs{5pt}

If $\Delta_{-1}=0$, by $(\ref{g1 -1 02})$, we must have $\Delta_1=3$ and $g_{1,-1}(\p,0)=a_{1,-1}^{0}$. Therefore, by $(\ref{g1 -1})$ we have $g_{1,-1}(\l,\p)=a_{1,-1}^{0}$, i.e., the first case of \eqref{eq-1} holds.\vs{5pt}

\noindent{\bf Case 2:} $\Delta_{-1}=1$.\vs{5pt}

In this case, $(\ref{g1 -1 02})$
shows that $\Delta_1=3$ and $a_{1,-1}^{1}=-\frac{a_{1,-1}^{0}}{\a_1}$. Using this in $(\ref{g1 -1})$, we get the second case of $(\ref{eq-1})$. \vs{5pt}

\noindent{\bf Case 3:} $\Delta_{-1}\neq0$ and $\Delta_{-1}\neq1$.
\vs{5pt}

If $\Delta_{-1}\neq0$ and $\Delta_{-1}\neq1$, comparing the degrees of $\p$ on both sides of $(\ref{g1 -1 02})$, we know that $g_{1,-1}(\p,0)=g_{1,-1}(0,0)\neq0$. Then comparing the coefficients of $\p^i$ for $i=1,2$ on both sides of $(\ref{g1 -1 02})$ respectively, we have $\Delta_{-1}+\Delta_{1}=3$ and $\Delta_{-1}+\Delta_{1}=1$, a contradiction.
Hence the lemma follows.
\hfill$\Box$ \vskip7pt

Now we deal with the case $\a_1 = 0$. If $a_{1,-1}^{0}\neq 0$, taking $\a_1 = 0$ in $(\ref{g1 -1 01})$ and $(\ref{g1 -1 02})$ and comparing the coefficients of $\p^2$ on both sides respectively, we can deduce that $\Delta_{-1}=0$, $\Delta_1=3$ and $g_{1,-1}(\p,0)=a_{1,-1}^{0}$.
Then by $(\ref{g1 -1})$ we immediately obtain the following.

\begin{lemm}\label{lemm3} If $\a_1 = 0$ and $a_{1,-1}^{0}\neq 0$, then $\Delta_{-1}=0$, $\Delta_1=3$ and $g_{1,-1}(\l,\p)=a_{1,-1}^{0}$.
\end{lemm}

Now we can consider the most complicated case that $\a_1 = 0$ and $a_{1,-1}^{0}= 0$.
In this case $(\ref{g1 -1 01})$ and $(\ref{g1 -1 02})$ turn into
\begin{eqnarray}\label{eq-2}
\!\!\!\!\!\!\!\!\!\!\!\!&\!\!\!\!\!\!\!\!\!\!\!\!\!\!\!&
(\Delta_{1}-1)(\Delta_{-1}-1)g_{1,-1}(\p,0)=0,
\\
\!\!\!\!\!\!\!\!\!\!\!\!&\!\!\!\!\!\!\!\!\!\!\!\!\!\!\!&\label{eq-3}
\Delta_{-1}(1-\Delta_{-1})g_{1,-1}(\p,0)=0.
\end{eqnarray}

\begin{lemm}\label{lemm4} If $\a_1 = 0$ and $a_{1,-1}^{0}= 0$, then we can deduce that $\Delta_{-1}=1,\ \Delta_1=3$ and
\begin{equation}\label{eq-4}
g_{1,-1}(\l,\p)=
a_{1,-1}^{1}(\l+\p)\ \ \ \ \  \mbox{for\ }\ \ \   a_{1,-1}^{1}\in\C^{*},
\end{equation}
\end{lemm}
\noindent{\it Proof.~}~By $(\ref{G_-1_j})$ and $(\ref{eq-3})$, our discussion will be divided into the following two cases.\vs{5pt}

\noindent{\bf Case 1:} $\Delta_{-1}=0$.\vs{5pt}

From $(\ref{eq-2})$, it follows that $(\Delta_{1}-1)g_{1,-1}(\p,0)=0$. If $\Delta_{1}\neq 1$, then it is obvious that $g_{1,-1}(\p,0)=0$. This together with $(\ref{g1 -1})$ shows that $g_{1,-1}(\l,\p)=0$, i.e., a contradiction with $(\ref{G_-1_j})$. Now we  suppose $\Delta_{1}=1$. Taking $\l=-\m\neq 0$ in $(\ref{g1 -1})$, using the fact that
$a_{1,-1}^{0}= 0$, we obtain
\begin{equation}\label{eq-5}
g_{1,-1}(\m,-\m)=-\frac{1}{2}\,g_{1,-1}(\m,0).
\end{equation}
In addition, $(\ref{g1 -1 0})$ shows that
\begin{equation}\label{eq-6}
g_{1,-1}(\m,\p)=g_{1,-1}(\m,0)\ \ \ \ \mbox{for\ \ } \p\neq -\m.
\end{equation}
Letting $\a_1=0$, $\Delta_{-1}=0$ and $\Delta_{1}=1$ in $(\ref{g1 -1 00})$, we have $\m g_{1,-1}(\m+\p,0)=(\m-2\p)g_{1,-1}(\m,0)$ for $\m+\p\neq 0$.
Then setting $\m=1$ and $\p=-\m+1$ in this formula respectively, we get
\begin{eqnarray}\label{eq-7}
\!\!\!\!\!\!\!\!\!\!\!\!&\!\!\!\!\!\!\!\!\!\!\!\!\!\!\!&
g_{1,-1}(\m,0)=(3-2\m)g_{1,-1}(1,0)\ \ \ \ \mbox{for\ \ } \m\neq 0,
\\
\!\!\!\!\!\!\!\!\!\!\!\!&\!\!\!\!\!\!\!\!\!\!\!\!\!\!\!&\label{eq-8}
\m g_{1,-1}(1,0)=(3\m-2)g_{1,-1}(\m,0).
\end{eqnarray}
Inserting $(\ref{eq-7})$ into $(\ref{eq-8})$ gives that $(\m-1)^{2}g_{1,-1}(1,0)=0$ for $\m\neq 0$. It leads to $g_{1,-1}(1,0)=0$. Then  $(\ref{eq-5})$, $(\ref{eq-6})$ together with $(\ref{eq-7})$ show that $g_{1,-1}(\l,\p)=0$, i.e., we also get a contradiction with $(\ref{G_-1_j})$.

\vs{5pt}
\noindent{\bf Case 2:} $\Delta_{-1}=1$.\vs{5pt}

Taking $\a_1=0$ and $\Delta_{-1}=1$ in $(\ref{g1 -1 00})$, we have
\begin{equation*}
\m\big((\Delta_{1}-1)\p-\m\big)g_{1,-1}(\m+\p,0)=(\m+\p)(2\p-\m)g_{1,-1}(\m,0).
\end{equation*}
Letting $\p=-\m+1$ in the above formula, then replacing $\m$ by $\p$, we can obtain
\begin{equation}\label{eq-9}
(2-3\p)g_{1,-1}(\p,0)=\p(\D_1-1-\D_1\p)g_{1,-1}(1,0).
\end{equation}
Recall that $g_{1,-1}(\p,0)=\sum_{i=1}^{m}a_{1,-1}^{i}\p^{i}$. Comparing the coefficients of $\p^i$ for $1\leq i\leq m+1$ on both sides of $(\ref{eq-9})$, we deduce that $a_{1,-1}^{i}=0$ for $2\leq i\leq m$, $a_{1,-1}^{1}=g_{1,-1}(1,0)$ and $(\D_1-3)a_{1,-1}^{1}=0$. Thus we have that
 $g_{1,-1}(\p,0)=a_{1,-1}^{1}\p$ and $(\D_1-3)a_{1,-1}^{1}=0$.
If $a_{1,-1}^{1}=0$, it follows that $g_{1,-1}(\p,0)=0$, then $(\ref{g1 -1})$ leads to $g_{1,-1}(\l,\p)=0$. Therefore, by $(\ref{G_-1_j})$ and $(\ref{g1 -1})$, we
have $g_{1,-1}(\l,\p)=a_{1,-1}^{1}(\l+\p)\neq 0$.
And this lemma holds.
\hfill$\Box$
\vskip12pt

\cl{\bf\S5. \  Classification of graded Lie conformal algebras}\setcounter{section}{5}\setcounter{equation}{0}\setcounter{theo}{0} \vs{5pt}

In this section,
we determine all $g_{i,j}(\l,\p)$ for $i,j\in\Z_{\geq-1}$, so that we can classify $\mathbb Z$-graded simple Lie conformal algebras ${\cal G}=\oplus_{i=-1}^\infty{\cal G}_i$.

\vs{5pt}


\noindent{\it Proof of Theorem \ref{main-theo}\,(3).~}~Applying ${\g_{-1}}_{\l}$ to $[\g{_0}_\m \g_{j}]$, by (\ref{J_a_b_c}), we have the Jacobi identity
$[\g{_{-1}}_{\l}[\g{_{0}}_\m \g_{j}]]=[[\g{_{-1}}_\l \g_{0}]_{\l+\m} \g_{j}]+[\g{_{0}}_\m[\g{_{-1}}_\l \g_j]]$.
Then comparing the coefficients of $\g_{j-1}$, by
$(\ref{G_i_j})$ and Lemma \ref{lemm1}, we obtain
\begin{eqnarray}
\label{eq-10'}
\!\!\!\!\!\!\!\!&\!\!\!\!\!\!\!\!&
(j\a_1+\l+\p+\D_j{\ssc\,}\m)g_{-1,j}(\l,\p)-\big((j-1)\a_1+\p+\D_{j-1}{\ssc\,}\m\big)g_{-1,j}(\l,\p+\m)
\nonumber\\\!\!\!\!\!\!\!\!\!\!\!\!&\!\!\!\!\!\!\!\!\!\!\!\!\!\!\!&\ \ \ \ \ \ \ \ \ \ \ \ \ \ \ \
=\big(\a_{1}+\l-(\D_{-1}-1)\m\big)g_{-1,j}(\l+\m,\p).
\end{eqnarray}
Taking $j=1$ in $(\ref{g-1 1 j})$, we have
\begin{eqnarray}
\label{eq-10}
\!\!\!\!\!\!\!\!&\!\!\!\!\!\!\!\!&
g_{-1,2}(\l,\p){\ssc\,}g_{1,1}(\m,\l+\p)-g_{-1,1}(\l,\p+\m){\ssc\,}g_{1,0}(\m,\p)
\nonumber\\\!\!\!\!\!\!\!\!\!\!\!\!&\!\!\!\!\!\!\!\!\!\!\!\!\!\!\!&\ \ \ \ \ \
=\big(\a_{1}+\p+\D_{1}{\ssc\,}(\l+\m)\big){\ssc\,}g_{-1,1}(\l,-\l-\m).
\end{eqnarray}
By skew-symmetry, we have $g_{-1,1}(\l,\p)=-g_{1,-1}(-\l-\p,\p)$, so for convenience, we suppose $a_{-1,1}^{0}=-a_{1,-1}^{0}$. By Lemma \ref{lemm2}--\ref{lemm4}, we need to consider the following three cases.
\vs{5pt}

\noindent{\bf Case 1:} $\Delta_{-1}=0$, $\Delta_{1}=3$ and $g_{1,-1}(\l,\p)=-a_{-1,1}^{0}\neq 0$.\vs{5pt}

Since we have skew-symmetry, by Lemma \ref{lemm1}, we have $g_{-1,1}(\l,\p)=a_{-1,1}^{0}$ and $g_{1,0}(\l,\p)=-\a_1+3\l+2\p$.
Then $(\ref{eq-10})$ turns into
\begin{equation*}
g_{-1,2}(\l,\p){\ssc\,}g_{1,1}(\m,\l+\p)=3a_{-1,1}^{0}(\l+2\m+\p).
\end{equation*}
Note that both $g_{-1,2}(\l,\p)$ and $g_{1,1}(\l,\p)$ are polynomials of $\l$ and $\p$, so by comparing the coefficients on both sides of the above formula, we can deduce that
\begin{eqnarray}
\label{eq-11}
\!\!\!\!\!\!\!\!&\!\!\!\!\!\!\!\!&
g_{-1,2}(\l,\p)=a_{-1,2}^{0},
\\\!\!\!\!\!\!\!\!\!\!\!\!&\!\!\!\!\!\!\!\!\!\!\!\!\!\!\!&
g_{1,1}(\l,\p)=\frac{3a_{-1,1}^{0}}{a_{-1,2}^{0}}(2\l+\p)\label{eq-12},
\end{eqnarray}
where $a_{-1,2}^{0}\in\C^*$.
Setting $j=2$ in $(\ref{eq-10'})$ and noting that $g_{-1,2}(\l,\p)=a_{-1,2}^{0}\neq 0$, we get
\begin{equation}
\label{eq-13}
\D_{2}=\D_{1}+1=4.
\end{equation}
By $(\ref{g-1 1 j})$, $(\ref{eq-10'})$ and $(\ref{eq-11})$--$(\ref{eq-13})$, we can inductively deduce that
\begin{eqnarray}
\label{eq-14}
\!\!\!\!\!\!\!\!&\!\!\!\!\!\!\!\!&
\D_{j}=j+2,
\\\!\!\!\!\!\!\!\!&\!\!\!\!\!\!\!\!&
g_{-1,j}(\l,\p)=a_{-1,j}^{0},\label{eq-15}
\\\!\!\!\!\!\!\!\!\!\!\!\!&\!\!\!\!\!\!\!\!\!\!\!\!\!\!\!&
g_{1,j}(\l,\p)=\frac{a_{-1,1}^{0}}{2a_{-1,j+1}^{0}}(j+2)\Big((j-1)\a_1+(j+3)\l+2\p\Big)\label{eq-16},
\end{eqnarray}
where $1\leq j\in\Z$ and $a_{-1,j}^{0}\in\C^*$.

Now we want to determine $g_{j,2}(\l,\p)$ for $2\leq j\in\Z$. Comparing the coefficients of $\g_{j+2}$ on both sides of
$[\g{_{j}}_{\l}[\g{_{1}}_\m \g_{1}]]=[[\g{_{j}}_\l \g_{1}]_{\l+\m} \g_{1}]+[\g{_{1}}_\m[\g{_{j}}_\l \g_1]]$, we obtain
\begin{eqnarray*}
\!\!\!\!\!\!\!\!&\!\!\!\!\!\!\!\!&
g_{j,2}(\l,\p){\ssc\,}g_{1,1}(\m,\l+\p)-g_{j,1}(\l,\p+\m){\ssc\,}g_{1,j+1}(\m,\p)
\nonumber\\\!\!\!\!\!\!\!\!\!\!\!\!&\!\!\!\!\!\!\!\!\!\!\!\!\!\!\!&\ \ \ \ \ \
=g_{j,1}(\l,-\l-\m){\ssc\,}g_{j+1,1}(\l+\m,\p).
\end{eqnarray*}From this,
using $(\ref{eq-16})$ and skew-symmetry, we obtain, for $1\leq j\in\Z$,
\begin{equation}
g_{j,2}(\l,\p)=\frac{-a_{-1,1}^{0}a_{-1,2}^{0}}{6a_{-1,j+1}^{0}a_{-1,j+2}^{0}}(j+2)(j+3)\Big((j-2)\a_1-(j+4)\l-(j+1)\p\Big)\label{eq-17}.
\end{equation}
Finally, we can determine all the polynomials $g_{j,i}(\l,\p)$ as follows. Noting from
$[\g{_{j}}_{\l}[\g{_{1}}_\m \g_{i}]]=[[\g{_{j}}_\l \g_{1}]_{\l+\m} \g_{i}]+[\g{_{1}}_\m[\g{_{j}}_\l \g_i]]$, we obtain
\begin{eqnarray}
\!\!\!\!\!\!\!\!&\!\!\!\!\!\!\!\!&
g_{j,i+1}(\l,\p){\ssc\,}g_{1,i}(\m,\l+\p)-g_{j,i}(\l,\p+\m){\ssc\,}g_{1,j+i}(\m,\p)
\nonumber\\\!\!\!\!\!\!\!\!\!\!\!\!&\!\!\!\!\!\!\!\!\!\!\!\!\!\!\!&\ \ \ \ \ \
=g_{j,1}(\l,-\l-\m){\ssc\,}g_{j+1,i}(\l+\m,\p).\label{eq-18}
\end{eqnarray}
By $(\ref{eq-16})$--$(\ref{eq-18})$, we can inductively deduce,  for $1\leq i\in\Z$, $1\leq j\in\Z$,
\begin{eqnarray}
\!\!\!\!\!\!\!\!&\!\!\!\!\!\!\!\!&
g_{j,i}(\l,\p)=\frac{-a_{-1,1}^{0}a_{-1,2}^{0}\cdots a_{-1,i}^{0}}{a_{-1,j+1}^{0}a_{-1,j+2}^{0}\cdots a_{-1,j+i}^{0}}
\times\frac{(j+2)(j+3)\cdots (j+i+1)}{(i+1)!}
\nonumber\\\!\!\!\!\!\!\!\!\!\!\!\!&\!\!\!\!\!\!\!\!\!\!\!\!\!\!\!&\ \ \ \ \ \ \ \ \ \ \ \ \
\ \ \times\Big((j-i)\a_1-(j+i+2)\l-(j+1)\p\Big).\label{eq-19}
\end{eqnarray}
In order for  the polynomials $g_{j,i}(\l,\p)$ to have some suitable forms, 
for $1\leq j\in\Z$, we replace $\g_j$ by $\g'_j=\frac{(j+1)!}{a_{-1,1}^{0}a_{-1,2}^{0}\cdots a_{-1,j}^{0}}\g_j$, so that $g_{-1,j}(\l,\p)$ and $g_{j,i}(\l,\p)$ have the following forms,
\begin{eqnarray*}
\!\!\!\!\!\!\!\!&\!\!\!\!\!\!\!\!&
g_{-1,j}(\l,\p)=j+1\mbox{ \ for \ }1\leq j\in\Z,
\\\!\!\!\!\!\!\!\!\!\!\!\!&\!\!\!\!\!\!\!\!\!\!\!\!\!\!\!&
g_{j,i}(\l,\p)=(i-j)\a_1+(i+j+2)\l+(j+1)\p\mbox{ \ for \ } 1\leq i\in\Z, 1\leq j\in\Z.
\end{eqnarray*}
Since in this case $\Delta_{-1}=0$, using $(\ref{g_0_j})$ and skew-symmetry, we obtain $g_{-1,0}(\l,\p)=\a_1-\p$.
By $(\ref{g_0_j})$ and $(\ref{eq-14})$, it is not hard to check that the second equation of the above also holds for $i=0$ or $j=0$.
Therefore, in this case we obtain $\GG\cong B(1,\a)$ for some $\a\in\C$. 
\vs{5pt}

\noindent{\bf Case 2:} $\Delta_{-1}=1$, $\Delta_{1}=3$, $\a_1=0$ and $g_{1,-1}(\l,\p)=a_{-1,1}^{1}(\l+\p)$, where $a_{-1,1}^{1}\in\C^*$.\vs{5pt}

In this case, using skew-symmetry, we obtain that $g_{-1,1}(\l,\p)=a_{-1,1}^{1}\l$ and $g_{1,0}(\l,\p)=3\l+2\p$.
Then $(\ref{eq-10})$ turns into
\begin{equation*}
g_{-1,2}(\l,\p){\ssc\,}g_{1,1}(\m,\l+\p)=3a_{-1,1}^{1}\l(\l+2\m+\p).
\end{equation*}
Since $g_{-1,2}(\l,\p)$ and $g_{1,1}(\l,\p)$ are polynomials of $\l$ and $\p$, we get from the above formula,
\begin{eqnarray}
\label{eq-20}
\!\!\!\!\!\!\!\!&\!\!\!\!\!\!\!\!&
g_{-1,2}(\l,\p)=a_{-1,2}^{1}\l,
\\\!\!\!\!\!\!\!\!\!\!\!\!&\!\!\!\!\!\!\!\!\!\!\!\!\!\!\!&
g_{1,1}(\l,\p)=\frac{3a_{-1,1}^{1}}{a_{-1,2}^{1}}(2\l+\p)\label{eq-21},
\end{eqnarray}
for some $a_{-1,2}^{1}\in\C^*$. Taking $j=2$, $\Delta_{-1}=1$, $\Delta_{1}=3$ and $\a_1=0$ in $(\ref{eq-10'})$, and noting that $a_{-1,2}^{1}\neq 0$, we get
\begin{equation}
\label{eq-22}
\D_{2}=4.
\end{equation}
Therefore, by $(\ref{g-1 1 j})$, $(\ref{eq-10'})$ and $(\ref{eq-20})$--$(\ref{eq-22})$, we can inductively prove that
\begin{eqnarray*}
\!\!\!\!\!\!\!\!&\!\!\!\!\!\!\!\!&
\D_{j}=j+2,
\\\!\!\!\!\!\!\!\!&\!\!\!\!\!\!\!\!&
g_{-1,j}(\l,\p)=a_{-1,j}^{1}\l,
\\\!\!\!\!\!\!\!\!\!\!\!\!&\!\!\!\!\!\!\!\!\!\!\!\!\!\!\!&
g_{1,j}(\l,\p)=\frac{a_{-1,1}^{1}}{2a_{-1,j+1}^{1}}(j+2)\Big((j+3)\l+2\p\Big),
\end{eqnarray*}
where $1\leq j\in\Z$ and $a_{-1,j}^{1}\in\C^*$.
Similar to Case 1, the above three formulae together with $(\ref{eq-18})$ inductively show the following,
\begin{eqnarray*}
\!\!\!\!\!\!\!\!&\!\!\!\!\!\!\!\!&
g_{j,i}(\l,\p)=\frac{a_{-1,1}^{1}a_{-1,2}^{1}\cdots a_{-1,i}^{1}}{a_{-1,j+1}^{1}a_{-1,j+2}^{1}\cdots a_{-1,j+i}^{1}}
\times\frac{(j+2)(j+3)\cdots (j+i+1)}{(i+1)!}
\nonumber\\\!\!\!\!\!\!\!\!\!\!\!\!&\!\!\!\!\!\!\!\!\!\!\!\!\!\!\!&\ \ \ \ \ \ \ \ \ \ \ \ \
\ \ \times\Big((j+i+2)\l+(j+1)\p\Big),
\end{eqnarray*}
for $1\leq i\in\Z$ and $1\leq j\in\Z$. Replace $\g_j$ by $\g'_j=\frac{(j+1)!}{a_{-1,1}^{1}a_{-1,2}^{1}\cdots a_{-1,j}^{1}}\g_j$ for $1\leq j\in\Z$, so that $g_{-1,j}(\l,\p)$ and $g_{j,i}(\l,\p)$ have the following forms,
\begin{eqnarray*}
\!\!\!\!\!\!\!\!&\!\!\!\!\!\!\!\!&
g_{-1,j}(\l,\p)=(j+1)\l\mbox{ \ for \ }1\leq j\in\Z,
\\\!\!\!\!\!\!\!\!\!\!\!\!&\!\!\!\!\!\!\!\!\!\!\!\!\!\!\!&
g_{j,i}(\l,\p)=(i+j+2)\l+(j+1)\p\mbox{ \ for \ } 1\leq i\in\Z,\, 1\leq j\in\Z.
\end{eqnarray*}
Using Lemma \ref{lemm1} and the fact that $\D_{j}=j+2$ for $1\leq j\in\Z$, we can immediately obtain that the above two formulae hold for all
$i,j\in\Z_{\geq-1}$. Therefore, we have proved that $g_{j,i}(\l,\p)=(i+j+2)\l+(j+1)\p$ for all $i,j\in\Z_{\geq-1}$, which is equivalent to that 
$\GG\cong B(2,0)$.
\vs{5pt}

\noindent{\bf Case 3:} $\Delta_{-1}=1$, $\Delta_{1}=3$, $\a_1\neq 0$ and $g_{1,-1}(\l,\p)=-\frac{a_{-1,1}^{0}}{\a_1}\big(\a_1-\l-\p\big)$, where $a_{-1,1}^{0}\in\C^*$.\vs{5pt}

In this case, by skew-symmetry, we get $g_{-1,1}(\l,\p)=a_{-1,1}^{0}\big(1+\frac{1}{\a_1}\l\big)$ and $g_{1,0}(\l,\p)=-\a_1+3\l+2\p$.
Then $(\ref{eq-10})$ leads to
\begin{equation*}
g_{-1,2}(\l,\p){\ssc\,}g_{1,1}(\m,\l+\p)=3a_{-1,1}^{0}(\l+2\m+\p)\Big(1+\frac{1}{\a_1}\l\Big).
\end{equation*}
Setting $\l=0$, $\p=0$ and $\l=\p=0$ in the above formula respectively, noting that
$g_{-1,2}(\l,\p)$ and $g_{1,1}(\l,\p)$ are polynomials of $\l$ and $\p$, we obtain
\begin{eqnarray}
\label{eq-23}
\!\!\!\!\!\!\!\!&\!\!\!\!\!\!\!\!&
g_{-1,2}(\l,\p)=a_{-1,2}^{0}\Big(1+\frac{1}{\a_1}\l\Big),
\\\!\!\!\!\!\!\!\!\!\!\!\!&\!\!\!\!\!\!\!\!\!\!\!\!\!\!\!&
g_{1,1}(\l,\p)=\frac{3a_{-1,1}^{0}}{a_{-1,2}^{0}}\big(2\l+\p\big)\label{eq-24},
\end{eqnarray}
where $a_{-1,2}^{0}\in\C^*$. Taking $j=2$, $\Delta_{-1}=1$ and $\Delta_{1}=3$ in $(\ref{eq-10'})$, noting that $a_{-1,2}^{0}\neq 0$, we get
\begin{equation}
\label{eq-25}
\D_{2}=4.
\end{equation}
By $(\ref{g-1 1 j})$, $(\ref{eq-10'})$ and $(\ref{eq-23})$--$(\ref{eq-25})$, we can inductively deduce
\begin{eqnarray*}
\!\!\!\!\!\!\!\!&\!\!\!\!\!\!\!\!&
\D_{j}=j+2,
\\\!\!\!\!\!\!\!\!&\!\!\!\!\!\!\!\!&
g_{-1,j}(\l,\p)=a_{-1,j}^{0}\Big(1+\frac{1}{\a_1}\l\Big),
\\\!\!\!\!\!\!\!\!\!\!\!\!&\!\!\!\!\!\!\!\!\!\!\!\!\!\!\!&
g_{1,j}(\l,\p)=\frac{a_{-1,1}^{0}}{2a_{-1,j+1}^{0}}(j+2)\Big((j-1)\a_1+(j+3)\l+2\p\Big),
\end{eqnarray*}
where $1\leq j\in\Z$ and $a_{-1,j}^{0}\in\C^*$.
Similar to Case 1, from $(\ref{eq-18})$ and the above three formulae, it inductively follows that
\begin{eqnarray*}
\!\!\!\!\!\!\!\!&\!\!\!\!\!\!\!\!&
g_{j,i}(\l,\p)=\frac{a_{-1,1}^{0}a_{-1,2}^{0}\cdots a_{-1,i}^{0}}{a_{-1,j+1}^{0}a_{-1,j+2}^{0}\cdots a_{-1,j+i}^{0}}
\times\frac{(j+2)(j+3)\cdots (j+i+1)}{(i+1)!}
\nonumber\\\!\!\!\!\!\!\!\!\!\!\!\!&\!\!\!\!\!\!\!\!\!\!\!\!\!\!\!&\ \ \ \ \ \ \ \ \ \ \ \ \
\ \ \times\Big((i-j)\a_1+(j+i+2)\l+(j+1)\p\Big),
\end{eqnarray*}
for $1\leq i\in\Z$ and $1\leq j\in\Z$. Since in this case $\a_1\neq 0$, by replacing $\g_j$ by $\g'_j=\frac{(j+1)!\a_1^{j}}{a_{-1,1}^{0}a_{-1,2}^{0}\cdots a_{-1,j}^{0}}\g_j$ for $1\leq j\in\Z$, we obtain that $g_{-1,j}(\l,\p)$ and $g_{j,i}(\l,\p)$ have the following forms,
\begin{eqnarray*}
\!\!\!\!\!\!\!\!&\!\!\!\!\!\!\!\!&
g_{-1,j}(\l,\p)=(j+1)(\a_1+\l)\mbox{ \ for \ } 1\leq j\in\Z,
\\\!\!\!\!\!\!\!\!\!\!\!\!&\!\!\!\!\!\!\!\!\!\!\!\!\!\!\!&
g_{j,i}(\l,\p)=(i-j)\a_1+(i+j+2)\l+(j+1)\p\mbox{ \ for \ } 1\leq i\in\Z,\, 1\leq j\in\Z.
\end{eqnarray*}
By Lemma \ref{lemm1} and noting that $\D_{j}=j+2$ for $1\leq j\in\Z$, we can immediately conclude that the above two formulae hold for all
$i,j\in\Z_{\geq-1}$. Hence, we obtain that $g_{j,i}(\l,\p)=(i-j)\a_1+(i+j+2)\l+(j+1)\p$ for all $i,j\in\Z_{\geq-1}$ and $\a_1\neq 0$. It follows that in this case, we have $\GG\cong B(2,\a)$ for some $\a\in\C^*$.

Therefore, the above three cases together show that Theorem \ref{main-theo} (3) holds.
\hfill$\Box$
\vs{10pt}


\end{CJK*}
\end{document}